\documentclass{amsart}

\usepackage{amsmath}
\usepackage{amssymb}
\usepackage{tikz}
\usepackage{subfig}
\usepackage{array}

\usetikzlibrary{calc}
\newtheorem{thm}{Theorem}[section]

\theoremstyle{remark}
\newtheorem{rem}[thm]{Remark}

\theoremstyle{definition}

\begin{document}

\title{Minimal Diagrams of Free Knots}
\author{Tomas Boothby}
\author{Allison Henrich}
\author{Alexander Leaf}
\address{Simon Fraser University, Department of Mathematics, Burnaby, B.C. V5A 1S6, Canada}
\address{Seattle University, Department of Mathematics, Seattle, WA 98122, United States}
\address{University of Michigan, Department of Mathematics, Ann Arbor, MI 48109, United States}
\date{\today}

\begin{abstract}
Manturov recently introduced the idea of a free knot, i.e. an equivalence class of virtual knots where equivalence is generated by crossing change and virtualization moves. He showed that if a free knot diagram is associated to a graph that is \emph{irreducibly odd}, then it is minimal with respect to the number of classical crossings. Not all minimal diagrams of free knots are associated to irreducibly odd graphs, however. We introduce a family of free knot diagrams that arise from certain  permutations that are minimal but not irreducibly odd. %We also explore the structure of the irreducibly odd graphs that give rise to minimal free knot diagrams in order to facilitate the classification of free knots.
\end{abstract}

\let\thefootnote\relax\footnotetext{{\em MSC classes}: 05C83, 57M27}

\maketitle

\section{Introduction \& Background} \label{sec:intro}
In recent years, many new theories of knots have become of particular interest. Louis Kauffman, for instance, introduced the notion of a virtual knot in his 1999 paper, \emph{Virtual Knot Theory}~\cite{VKT}. Virtual knot theory is a generalization of knot theory that has several interpretations: virtual knots can be thought of as equivalence classes of knot diagrams that may contain a new type of ``virtual" crossing; they can be viewed as knots in thickened surfaces (where classical knots are viewed as knots in thickened spheres); and they can be thought of as equivalence classes of Gauss diagrams or Gauss codes. 

Viewing virtual knots as Gauss diagrams, several other knot theories have come to life, among them the theories of virtual strings and free knots. Virtual strings were introduced in 2004 by Turaev~\cite{Turaev}, and the similar notion of free knots was introduced in 2009 by Manturov~\cite{FreeKnots},~\cite{FreeKnots2}. Both generalizations come from first viewing virtual knots as Gauss diagrams, then ``forgetting" certain information contained in the diagram. This forgetting produces new equivalences between previously distinct virtual knots. Nonetheless, highly non-trivial theories of knots remain.

\begin{figure}[ht]
\begin{tikzpicture}[thick,scale=1.5]
 \draw (0,0) to [out=0,in=180] (1,.5);
 \draw (0,.5) to [out=0,in=180] (1,0);

 \draw (2,0)  to [out=0,in=180]
       (2+1/3,.5) to [out=0,in=180]
       (2+2/3,0) to [out=0,in=180]
       (3,.5);
 \draw (2,.5) to [out=0,in=180] 
       (2+1/3,0)  to [out=0,in=180]
       (2+2/3,.5)  to [out=0,in=180]  
       (3,0);

 \draw [thin] (2+1/6,.25) circle (.0625/1.5);
 \draw [thin] (3-1/6,.25) circle (.0625/1.5);
 \node (arrow) at (1.5,.25) {$\leftrightarrow$};
\end{tikzpicture}
\caption{The flat virtualization move}
\label{fig:flat_vm}
\end{figure}
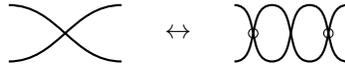

In this paper, we focus our attention on the theory of free knots. Free knots can be defined as equivalence classes of knot diagrams with two types of crossings: flat classical crossings and virtual crossings (denoted by an encircled flat crossing). Two diagrams of free knots are equivalent if they can be related by a sequence of flat classical and virtual Reidemeister moves as well as the flat virtualization move. See Figures~\ref{fig:flat_rm}~and~\ref{fig:flat_vm}. 

In a sense, flat classical crossings are akin to ordinary classical crossings where we have ``forgotten'' which strand is on top. The virtualization move is necessary if we ``forget'' from which direction we ought to approach a given classical crossing in our knot diagram.

\begin{figure}[ht]
\begin{tabular}{ccc}
\begin{tikzpicture}[thick,scale=1.5]
 \draw (-1,1) to (-1,.9);
 \draw (-1,.1) to (-1,0);
 \draw (-1,.9)   [out=270,in=45] 
    to (-1.15,.5) [out=225,in=270] 
    to (-1.45,.5)  [out=90, in=135] 
    to (-1.15,.5) [out=315,in=90] 
    to (-1,.1);
 \node (arrow) at (-.5,.5) {$\leftrightarrow$};
 \draw [thin] (-1.15,.5) circle (.0625/1.5);

 \draw (0,1) to (0,0);
 \draw (1,1) to (1,.9);
 \draw (1,.1) to (1,0);
 \draw (1,.9)   [out=270,in=135] 
    to (1.15,.5) [out=315,in=270] 
    to (1.45,.5)  [out=90, in=45] 
    to (1.15,.5) [out=225,in=90] 
    to (1,.1);
 \node (arrow) at (.5,.5) {$\leftrightarrow$};
\end{tikzpicture}
&\quad&
\begin{tikzpicture}[thick,scale=1.5]
 \draw (0,1) to (0,0);
 \draw (-.2,1) to (-.2,0);

 \draw (1,1) to [out=270,in=90]
       (1.2,.5) to [out=270,in=90]
       (1,0);
 \draw (1.2,1) to [out=270,in=90]
       (1,.5) to [out=270,in=90]
       (1.2,0);
 
 \draw (-1,1) to [out=270,in=90]
       (-1.2,.5) to [out=270,in=90]
       (-1,0);
 \draw (-1.2,1) to [out=270,in=90]
       (-1,.5) to [out=270,in=90]
       (-1.2,0);

 \draw [thin] (-1.1,.25) circle (.0625/1.5); 
 \draw [thin] (-1.1,.75) circle (.0625/1.5);

 \node (arrow) at (-.5,.5) {$\leftrightarrow$};

 \node (arrow) at (.5,.5) {$\leftrightarrow$};
\end{tikzpicture}
\end{tabular}

\begin{tabular}{ccccc}
\begin{tikzpicture}[thick,scale=1.5]
 \draw (0,  1) to [out=270,in=90] (-.4,0);
 \draw (-.2,1) to [out=270,in=90] (-.4,.5) to [out=270,in=90] (-.2,0);
 \draw (-.4,1) to [out=270,in=90] (0,  0);

 \draw (1,  1) to [out=270,in=90] (1.4,0);
 \draw (1.2,1) to [out=270,in=90] (1.4,.5) to [out=270,in=90] (1.2,0);
 \draw (1.4,1) to [out=270,in=90] (1,  0);
 \node (arrow) at (.5,.5) {$\leftrightarrow$};
\end{tikzpicture}
&\quad&
\begin{tikzpicture}[thick,scale=1.5]
 \draw (0,  1) to [out=270,in=90] (-.4,0);
 \draw (-.2,1) to [out=270,in=90] (-.4,.5) to [out=270,in=90] (-.2,0);
 \draw (-.4,1) to [out=270,in=90] (0,  0);

 \draw (1,  1) to [out=270,in=90] (1.4,0);
 \draw (1.2,1) to [out=270,in=90] (1.4,.5) to [out=270,in=90] (1.2,0);
 \draw (1.4,1) to [out=270,in=90] (1,  0);
 
 \draw [thin] (-.2,.5) circle (.0625/1.5);
 \draw [thin] (1.2,.5) circle (.0625/1.5);

 \draw [thin] (1.325,.7125) circle (.0625/1.5);
 \draw [thin] (1.325,.2875) circle (.0625/1.5);

 \draw [thin] (-.325,.7125) circle (.0625/1.5);
 \draw [thin] (-.325,.2875) circle (.0625/1.5);

 \node (arrow) at (.5,.5) {$\leftrightarrow$};
\end{tikzpicture}
&\quad&
\begin{tikzpicture}[thick,scale=1.5]
 \draw (0,  1) to [out=270,in=90] (-.4,0);
 \draw (-.2,1) to [out=270,in=90] (-.4,.5) to [out=270,in=90] (-.2,0);
 \draw (-.4,1) to [out=270,in=90] (0,  0);

 \draw (1,  1) to [out=270,in=90] (1.4,0);
 \draw (1.2,1) to [out=270,in=90] (1.4,.5) to [out=270,in=90] (1.2,0);
 \draw (1.4,1) to [out=270,in=90] (1,  0);
 
 %\draw (-.2,.5) circle (.0625);
 %\draw (1.2,.5) circle (.0625);

 \draw [thin] (1.325,.7125) circle (.0625/1.5);
 \draw [thin] (1.325,.2875) circle (.0625/1.5);

 \draw [thin] (-.325,.7125) circle (.0625/1.5);
 \draw [thin] (-.325,.2875) circle (.0625/1.5);

 \node (arrow) at (.5,.5) {$\leftrightarrow$};
\end{tikzpicture}
\end{tabular}

\caption{Flat classical and virtual Reidemeister moves}
\label{fig:flat_rm}
\end{figure}
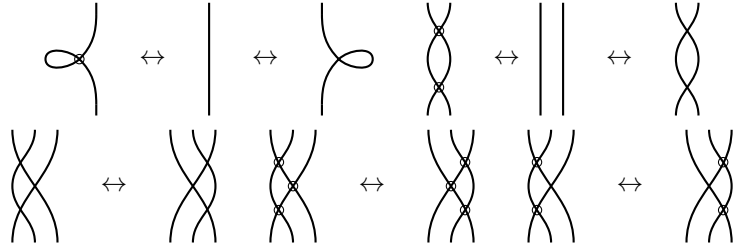

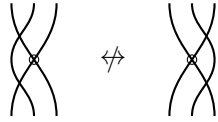
\begin{figure}[ht]
\begin{tikzpicture}[thick,scale=1.5]
 \draw (0,  1) to [out=270,in=90] (-.4,0);
 \draw (-.2,1) to [out=270,in=90] (-.4,.5) to [out=270,in=90] (-.2,0);
 \draw (-.4,1) to [out=270,in=90] (0,  0);

 \draw (1,  1) to [out=270,in=90] (1.4,0);
 \draw (1.2,1) to [out=270,in=90] (1.4,.5) to [out=270,in=90] (1.2,0);
 \draw (1.4,1) to [out=270,in=90] (1,  0);
 
 \draw [thin] (-.2,.5) circle (.0625/1.5);
 \draw [thin] (1.2,.5) circle (.0625/1.5);
%D73645935
 \node (arrow) at (.5,.5) {$\leftrightarrow$};
 \node (not) at (.5,.5) {$/$};
\end{tikzpicture}
\caption{The forbidden move}\label{forbidden}
\end{figure}

Note that the move pictured in Figure~\ref{forbidden} is forbidden for free knots. We also note that all free knots with either no virtual crossings or no classical crossings are trivial.

\begin{figure}
  \begin{center}
    \begin{tabular}{m{1.2in}m{1.2in}m{1.2in}}
      \begin{tikzpicture}[scale=1]
          \coordinate (0) at (90:.5);
          \coordinate (1) at (210:.5);
          \coordinate (2) at (330:.5);
          \coordinate (3) at (90:1);
          \coordinate (4) at (210:1);
          \coordinate (5) at (330:1);
          \draw[thick] (0) -- (1);
          \draw[thick] (0) -- (2);
          \draw[thick] (0) -- (3);
          \draw[thick] (1) -- (2);
          \draw[thick] (1) -- (4);
          \draw[thick] (2) -- (5);
          \fill (0) circle (1.5pt);
          \fill (1) circle (1.5pt);
          \fill (2) circle (1.5pt);
          \fill (3) circle (1.5pt);
          \fill (4) circle (1.5pt);
          \fill (5) circle (1.5pt);
      \end{tikzpicture}
      &
    \begin{tikzpicture}
      \foreach \x in {0, ..., 11} \coordinate (\x) at (30.0000*\x:1);
      \draw [very thick] (0,0) circle (1);
      \draw [thick]
            (0) [out=180.0000,in=480.0000] to (10)
            (1) [out=210.0000,in=390.0000] to (7)
            (2) [out=240.0000,in=300.0000] to (4)
            (3) [out=270.0000,in=450.0000] to (9)
            (5) [out=330.0000,in=510.0000] to (11)
            (6) [out=360.0000,in=420.0000] to (8)
        ;
    \end{tikzpicture}
      &
\begin{tikzpicture}[scale=.5]
 
  \coordinate (6) at (0,0);
  \coordinate (7) at (-1,1.5);
  \coordinate (4) at (1,1.5);
  \coordinate (0) at (-2,2);
  \coordinate (3) at (2,2);
  \coordinate (8) at (-1,2.5);
  \coordinate (5) at (1,2.5);
  \coordinate (9) at (0,3.5);
  \coordinate (1) at (-.5,4.25);
  \coordinate (2) at (.5,4.25);

  \draw (2) circle (.125);
  \draw (4) circle (.125);
  \draw (6) circle (.125);
  \draw (8) circle (.125);

  \draw[thick] (1)  [out=0,in=180] to
               (2)  [out=0,in=45] to
               (3)  [out=225,in=0] to
               (4)  [out=180,in=180] to
               (5)  [out=0,in=135] to
               (3)  [out=315,in=315] to
               (6)  [out=135,in=270] to
               (7)  [out=90,in=270] to
               (8)  [out=90,in=225] to
               (9)  [out=45,in=270] to
               (2)  [out=90,in=90] to
               (1)  [out=270,in=135] to
               (9)  [out=315,in=90] to
               (5)  [out=270,in=90] to
               (4)  [out=270,in=45] to
               (6)  [out=225,in=225] to
               (0)  [out=45,in=180] to
               (8)  [out=0,in=0] to
               (7)  [out=180,in=315] to
               (0)  [out=135,in=180] to
               (1);

\end{tikzpicture}
    \end{tabular}
  \end{center}
  \caption{The 3-morningstar, its chord diagram, and its free knot diagram}\label{triskelion}
\end{figure}
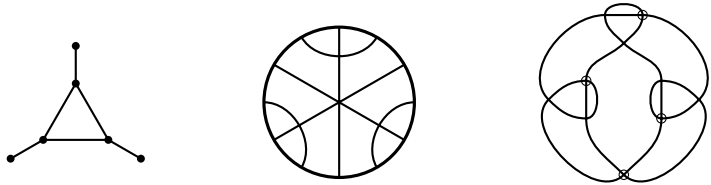

It was shown in~\cite{FreeKnots} that if a certain graph (called the {\em intersection graph}) associated to a free knot diagram is {\em irreducibly odd}, then the corresponding diagram is {\em minimal} in the sense that it has the fewest number of classical crossings of all free knots in its equivalence class. Manturov defines a graph to be \emph{odd} if every vertex $v \in G$ has odd degree and
\emph{irreducibly odd} if, in addition, for any pair of vertices $u,v \in G$ there is a third vertex $w$ adjacent
to either $u$ or $v$ but not both. Figure~\ref{triskelion} illustrates the simplest example of a free knot corresponding to an irreducibly odd graph, the 3-morningstar. We note that this is the only example of a free knot with fewer than seven classical crossings that is associated to an irreducibly odd graph. 

To find the intersection graph associated to a free knot diagram, it is useful to first draw the associated chord diagram. Given a free knot, we draw its chord diagram by numbering the classical crossings in the diagram. We then, separately, draw a parametrizing (core) circle. We choose a base point on the knot and a corresponding base point on the circle. We then traverse the knot (in the direction of a chosen orientation) from the base point, keeping track along the circle of the order in which we encounter the crossings. Each crossing will appear twice along the circle, so we are able to form one chord in the circle corresponding to each crossing. For those familiar with Gauss diagrams, we note that the chord diagram can be derived from the Gauss diagram simply by omitting the arrowhead and the sign on each chord.

Once we have a chord diagram, we form our corresponding intersection graph. For every chord in the chord diagram, we have a vertex in our graph. There is an edge between two vertices if and only if the two corresponding chords in the chord diagram intersect. 

\begin{figure}[h]
    \begin{tikzpicture}[font=\small]
        \node (a) at (0.0000,3.1235) {c};
        \node (b) at (1.2494,2.4988) {e};
        \node (c) at (0.0000,1.8741) {b};
        \node (d) at (1.2494,1.2494) {d};
        \node (e) at (2.4988,1.8741) {g};
        \node (f) at (2.4988,0.6247) {h};
        \node (g) at (2.4988,3.1235) {f};
        \node (h) at (0.0000,0.6247) {a};
        \draw[thick] (h) to (c) to (a) to (b) to (g) to (e) to (d) to (c);
        \draw[thick] (a) to (g);
        \draw[thick] (e) to (f);
        \draw[thick] (b) to (d);
    \end{tikzpicture}
    \begin{tikzpicture}[scale=1.3,thick,font=\tiny]
        \draw [very thick] (0,0) circle (1);
        \foreach \x in {0,1,...,15}
            \coordinate (\x) at (\x * 22.5:1);
        \draw (13) [out=180+13*22.5,in=180+15*22.5]  to (15);  \node (a) at (15*22.5:1.1)  {a};
        \draw (2) [out=180+2*22.5,in=180+14*22.5] to (14); \node (b) at (2*22.5:1.1)  {b};
        \draw (3) [out=180+3*22.5,in=180+7*22.5]  to (7);  \node (f) at (7*22.5:1.1)  {f};
        \draw (4) [out=180+4*22.5,in=180+12*22.5] to (12); \node (e) at (12*22.5:1.1) {e};
        \draw (5) [out=180+5*22.5,in=180+1*22.5] to (1); \node (c) at (1*22.5:1.1) {c};
        \draw (6) [out=180+6*22.5,in=180+10*22.5] to (10); \node (g) at (6*22.5:1.1)  {g};
        \draw (8) [out=180+8*22.5,in=180+0*22.5] to (0); \node (d) at (0*22.5:1.1) {d};
        \draw (9) [out=180+9*22.5,in=180+11*22.5] to (11); \node (h) at (9*22.5:1.1)  {h};
    \end{tikzpicture}
\begin{tikzpicture}[scale=.5,font=\tiny]

  \coordinate (0) at (2.5, .5);
  \coordinate (1) at (3.5, 3 + 1/3);
  \coordinate (2) at (3,   4);
  \coordinate (3) at (5,   5);
  \coordinate (4) at (1.5, 3);
  \coordinate (5) at (0,   5);
  \coordinate (a) at (4,   4);
  \coordinate (b) at (5,   4);
  \coordinate (c) at (2.5, 1.5);
  \coordinate (d) at (3.5, 2 + 2/3);
  \coordinate (e) at (2,   4);
  \coordinate (f) at (2.5, 5);
  \coordinate (g) at (0,   4);
  \coordinate (h) at (1,   4);

  \draw (0) circle (.125);
  \draw (1) circle (.125);
  \draw (2) circle (.125);
  \draw (3) circle (.125);
  \draw (4) circle (.125);
  \draw (5) circle (.125);

  \draw[thick](g) [out=0,in=180] to 
              (h) [out=0,in=180] to 
              (e) [out=0,in=180] to 
              (2) [out=0,in=180] to 
              (a) [out=0,in=180] to 
              (b) [out=0,in=0] to 
              (3) [out=180,in=90] to 
              (a) [out=270,in=45] to 
              (1) [out=225,in=135] to 
              (d) [out=315,in=45] to 
              (c) [out=225,in=135] to 
              (0) [out=315,in=270] to 
              (b) [out=90,in=270] to 
              (3) [out=90,in=45] to 
              (f) [out=225,in=90] to 
              (e) [out=270,in=45] to 
              (4) [out=225,in=135] to 
              (c) [out=315,in=45] to 
              (0) [out=225,in=270] to 
              (g) [out=90,in=270] to 
              (5) [out=90,in=135] to 
              (f) [out=315,in=90] to 
              (2) [out=270,in=135] to 
              (1) [out=315,in=45] to 
              (d) [out=225,in=315] to 
              (4) [out=135,in=270] to 
              (h) [out=90,in=0] to 
              (5) [out=180,in=180] to 
              (g);

  \node (A) at (4.2, 4.2) {a};
  \node (B) at (5.2, 3.8) {b};
  \node (C) at (2.5, 1.8) {c};
  \node (D) at (3.8, 2 + 2/3) {d};
  \node (E) at (2.2,   3.8) {e};
  \node (F) at (2.5, 5.3) {f};
  \node (G) at (-.2,   3.8) {g};
  \node (H) at (.8,   3.8) {h};
\end{tikzpicture}
    \caption{An example of a free knot with its corresponding chord diagram and graph}\label{ex}
\end{figure}
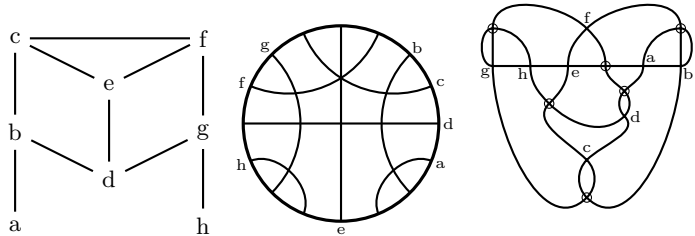

In Figure~\ref{ex}, we provide another example of an irreducibly odd graph and its corresponding free knot diagram and chord diagram. %For many more examples of minimal diagram free knots, we refer the reader to the irreducibly odd graphs in Table~\ref{table} that are associated to chord diagrams.

Clearly, every free knot has an associated graph. We note, however, that there are many graphs in general and irreducibly odd graphs in particular that do not have associated chord diagrams, and hence, cannot be associated to free knots.  

In this paper, we define the notion of a permutation graph. We use this notion to provide examples of minimal free knot diagrams whose graphs are {\em not} irreducibly odd.

\section{Permutation Graphs}\label{permutation}

If a free knot diagram is associated to an irreducibly odd graph, then it is a minimal crossing representative of the free knot it represents. The converse of this statement fails to be true. Indeed, Manturov gives a counterexample to the converse in~\cite{FreeKnots2}. Here, we generalize Manturov's example to a family of minimal free knot diagrams that are associated with permutations. Before we can describe this family, we introduce several definitions.

\subsection{Definitions}

A \emph{permutation chord diagram} is a chord diagram such that it is possible to add a chord to the diagram that intersects every other chord. This chord is called the \emph{equator} and is often denoted by a dashed line, when pictured. A \emph{permutation graph} is a graph that is associated to a permutation chord diagram.  Throughout this paper, we will write permutations in one-line notation, where $[\sigma_1 ~ \sigma_2 ~ \cdots ~ \sigma_n]$ denotes the function $\sigma(i)=\sigma_i$ for $i=1, 2, \cdots, n$.

Following the terminology from the study of permutations, we say that $\{i,j\}$ is an \emph{inversion} if the chords corresponding to $i$ and $j$ intersect in the chord diagram. This is the case if and only if the two corresponding vertices in the permutation graph are adjacent.

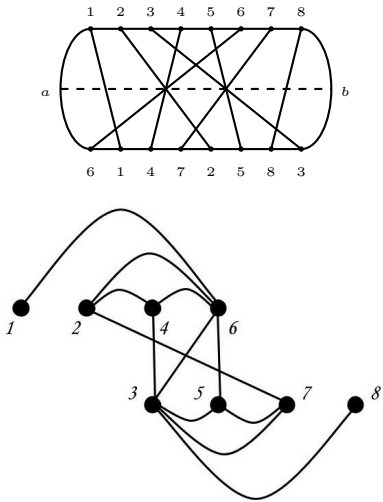
\begin{figure}
\begin{center}   \begin{center}
   \newcommand{\chordheight}{2}
   \begin{tikzpicture}[thick,text height=1.5ex,text depth=.25ex,x=.4cm,y=.8cm]
    \foreach \i/\j in {1/2,2/5,3/8,4/3,5/6,6/1,7/4,8/7} {
      \coordinate[label=above:{$\scriptscriptstyle \i$}] (a\i) at (\i,\chordheight);
      \coordinate[label=below:{$\scriptscriptstyle \i$}] (b\i) at (\j,0);
      \draw (a\i) to (b\i);
    }
    \coordinate[label=right:{$\scriptscriptstyle b$}] (b) at (8+\chordheight/2,\chordheight/2);
    \coordinate[label=left:{$\scriptscriptstyle a$}] (a) at (1-\chordheight/2,\chordheight/2);
    \draw[dashed] (a) to (b);
    \draw (1,0) arc (270:90:\chordheight/2);
    \draw (8,\chordheight) arc (90:-90:\chordheight/2);
    \draw (1,0) to (8,0);
    \draw (1,\chordheight) to (8,\chordheight);

    \foreach \i in {1,...,8} {
      \fill (\i,0) circle (1pt); 
      \fill (\i,\chordheight) circle (1pt);
     } 
   \end{tikzpicture}
  \end{center}\includegraphics[height=1.7in]{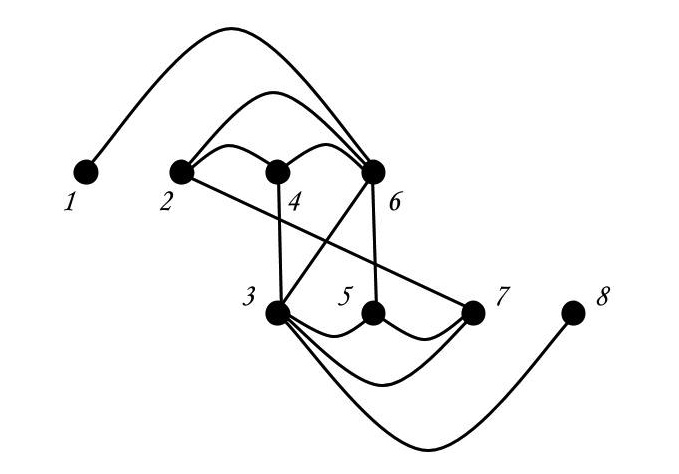}\end{center}
\caption{The chord diagram and intersection graph for $\sigma = [61472583]$.}\label{square}
\end{figure}

Note that the 3-morningstar does \emph{not} correspond to a permutation chord diagram, while the chord diagram in Figure~\ref{square} is a permutation chord diagram.

%\begin{figure}[h]
%\begin{center}\includegraphics[height=1.5in]{perm_ex.png}\end{center}
%\caption{The permutation chord diagram for $[3,4,1,2]$ and its associated permutation graph}\label{square}
%\end{figure}

To see how a permutation chord diagram with $n$ chords corresponds to an element (or, more precisely, elements) of the symmetric group $S_n$, we draw an equator in our chord diagram. We typically rotate the circle so that the equator is horizontal with left endpoint labelled $a$ and right endpoint labelled $b$, as in Figure~\ref{square}. We then traverse the circle from $a$ to  $b$ along the northern hemisphere, labeling endpoints of chords $1, 2, \cdots, n$ as we encounter them. Once we have assigned labels to each chord, we traverse the circle from $a$ to $b$ in the opposite direction (along the southern hemisphere), noting in which order we encounter our numbered chords. This reordering of the numbers $1, 2,\cdots, n$ gives a permutation, $\sigma$. To be precise, $\sigma(i)$ is the label of the $i$th chord we encounter when traversing the circle from $a$ to $b$ along the southern hemisphere of the chord diagram.

Note that if we started by traversing the circle from $a$ to $b$ in the opposite direction, we would obtain the inverse permutation, $\sigma ^{-1}$. If we instead traversed the circle from $b$ to $a$, the permutation we would get is the conjugate $\alpha \sigma \alpha ^{-1}$ of $\sigma$, where $\alpha$ is the order-reversing permutation $[n, n-1, \cdots, 2, 1]$.  Finally, we note that some chord diagrams may have multiple distinct equators, possibly yielding distinct permutations. %In fact, Figure~\ref{square} is such an example.

From this description, it is easy to see how to construct the permutation chord diagram associated to a given permutation $\sigma$. While we have seen that there are several related yet distinct permutations that may correspond to a given chord diagram, we remark that a permutation defines a unique chord diagram.

\subsection{Permutations and Minimal Free Knot Diagrams}

It is natural now to consider the relationship between irreducibly odd graphs and permutation graphs. To see how the two notions are related, we discuss several properties of permutation graphs related to irreducible oddness.

We say that a permutation $\sigma$ on $\{1, \cdots, n\}$ is \emph{parity-preserving} if $\sigma_{i}$ is even if and only if $i$ is even. Moreover, we call $\sigma$ \emph{parity-reversing} if $\sigma_{i}$ is even if and only if $i$ is odd. Note that a permutation graph that is irreducibly odd must be parity-reversing because an irreducibly odd permutation must be odd.  

Furthermore, a permutation $\sigma$ that yields an irreducibly odd graph must be \emph{nonconsecutive}: that is, for all $i$, $| \sigma_{i + 1} - \sigma_i |>1$. When we view the permutation as a chord diagram, we see that if $| \sigma_{i + 1} - \sigma_i |=1$, then for any chord $j\neq i, i+1$, the chords $i$ and $j$ intersect if and only if $i+1$ and $j$ intersect. This violates irreducibility, so no such adjacent pair can exist.

While the parity-reversing and nonconsecutive properties are indeed necessary for a permutation graph to be irreducibly odd, they are not sufficient.  For example, consider the permutation $\sigma = [10, 3, 6, 9, 4, 7, 2, 5, 8, 1]$.  Observe that $\sigma$ is nonconsecutive and parity-reversing, but fails to be irreducibly odd. This is because $\{10, i\}$ and $\{1, i\}$ are both inversions for $i \not= 1, 10$. This means that the vertices corresponding to 1 and 10 in the permutation graph are each adjacent to every other vertex. Hence, irreducibility is not satisfied.

Finally, we say that a permutation $\sigma$ on $\{ 1, \cdots, n \}$ is \emph{crossed} if there exists a unique $i$ such that $\{ j , i \}$ is an inversion for all $j \not= i$. Note that if $\sigma$ is crossed with crossing element 1, then $\sigma_n = 1$.  Furthermore, $ \tau = [ \sigma_1, \cdots, \sigma_{n-1} ]$ is a permutation on $\{ 2, \cdots , n \} $ that is not crossed.  

It is clear that any permutation chord diagram that corresponds to a crossed permutation $\sigma$ can be relabeled as a chord diagram where the crossing element is 1. We simply find the crossing element $i$ and relabel it as chord $1$, and moving clockwise around the circle, label the remaining chords $2, \cdots, n$.

\begin{rem}Note that, for non-trivial free knots with crossed permutation chord diagrams, there are no available simplifying type 1 Reidemeister moves and any available classical type 3 moves would involve the crossing element.\end{rem}

\begin{thm}\label{min}
Suppose that $\sigma$ is a crossed nonconsecutive permutation on $\{ 1, 2, \cdots, n\}$ with crossing element 1 satisfying 
\begin{enumerate}
\item $|\sigma _{i-1}-\sigma _{i}|<n-2$ for $i=2, \cdots, n-1$, and
\item $1<|\sigma _1-\sigma _{n-1}|<n-2$.
\end{enumerate}
Then $\sigma$ corresponds to a minimal free knot diagram.
\end{thm}

For ease of exposition, let us refer to a permutation that satisfies the conditions of Theorem~\ref{min} as a {\em minimal permutation}.

The reader may verify that Manturov's example in~\cite{FreeKnots2} corresponds to the permutation $[5,2,7,4,9,6,3,8,1]$, which satisfies the conditions of the theorem. We see that this permutation is parity-preserving, so it fails to be irreducibly odd.

Before we begin the proof of Theorem~\ref{min}, we recall several useful free knot and link invariants from Manturov in~\cite{FreeKnots},~\cite{FreeKnots2} that were inspired by work of Turaev and Goldman.

\begin{figure}[h]
\begin{center}\includegraphics[height=1in]{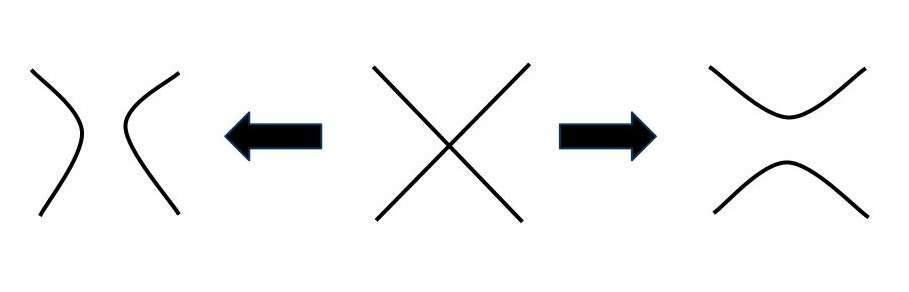}\end{center}
\caption{Smoothing at a crossing: one smoothing produces a knot and the other a 2-component link.}\label{smooth}
\end{figure}

Suppose $G$ is a fixed diagram of a free knot. Let $\Delta$ denote the sum $$\Delta(G)=\sum_c[G_c].$$ Here, $c$ denotes a crossing in the free knot diagram and $G_c$ is the the free two-component link diagram obtained by performing the smoothing on $c$ that results in a link. See Figure~\ref{smooth} for an illustration of smoothing at a crossing. The equivalence relation $\sim$ that defines the equivalence class $[\cdot]$ is generated by the flat Reidemeister 2 move only, with all diagrams containing free loops taken to be 0. Finally, we consider $\Delta (G)$ to be a formal sum of equivalence classes of free two-component link diagrams with $\mathbb{Z}_2$ coefficients. The resulting value of $\Delta(G)$ is independent of the particular diagram we chose for $G$, so $\Delta$ is an invariant of free knots.

Now, we recall Manturov's invariant on two-component free links. Denote by the bracket $\{ \cdot\}$ the invariant given by $$\{L\}=\sum_{s}L_s.$$ Here, $L_s$ is the free knot or link diagram obtained by smoothing at each crossing involving a single component of $L$. Once again, equivalence is generated by the flat Reidemeister type 2 move only. All possible combinations of smoothing such crossings (both with and against orientation) should be included in the sum, however, any $L_s$ containing a simple free loop is taken to be 0 and coefficients are in $\mathbb{Z}_2$.

%Add some pictures to make these invariants clear
Armed with these invariants, we proceed with our proof. We follow the strategy used in the proof of Statement 2 in~\cite{FreeKnots2}. Statement 2 shows that a particular example of a crossed permutation that satisfies the required properties corresponds to a minimal free knot diagram. Here, we generalize the result.

\begin{proof}[Proof of Theorem~\ref{min}]
Let $\sigma$ be a permutation as described in the hypotheses of the theorem. Let $G$ be the corresponding free knot diagram. We consider the invariant $$\Delta (G)=[G_1]+[G_2]+\cdots [G_n]$$ where $G_i$ refers to the free link obtained by smoothing at crossing $i$.
Because $\{1,j\}$ is an inversion for all $j\neq 1$ (since $1$ is the crossing element of the crossed permutation $\sigma$), each $j$ represents a crossing in $G_1$ that involves \emph{both} components in the two-component free link. Thus, $\{G_1\}=G_1$. 

Note that $G_1$ is a minimal diagram in the equivalence class generated by the flat Reidemeister 2 move. Indeed, the properties of $\sigma$ guarantee that no Reidemeister 2 moves are available after the smoothing. Property (1) ensures that chords $2$ and $n$ cannot be removed by a type 2 move after the smoothing, while Property (2) guarantees that chords $\sigma_1$ and $\sigma_{n-1}$ cannot be removed by a type 2 move after the smoothing. The nonconsecutive property guarantees that no other pairs can be removed. Thus, no diagram that is equivalent to $G_1$ by a sequence of type 2 moves has fewer than $n-1$ crossings.

On the other hand, the number of crossings involving a single component in any $G_j$ with $j\neq 1$ is positive since 1 is the only crossing element in $\sigma$. So $\{G_j\}$ does not have $G_1$ as one of its terms. This follows from the fact that each link in $\{G_j\}$ has a diagram with strictly fewer than $n-1$ crossings. Thus, $\{\Delta (G)\}$ contains the term $G_1$. Since this term has $n-1$ crossings and $\{\Delta (G)\}$ is an invariant, the free knot $G$ must have a minimum of $n$ crossings. Moreover, it was shown by Manturov that, up to purely virtual Reidemeister moves, there is a unique free knot diagram that realizes the minimum crossing number.

\end{proof}

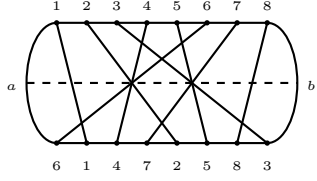
\begin{figure}
  \begin{center}
   \newcommand{\chordheight}{2}
   \begin{tikzpicture}[thick,text height=1.5ex,text depth=.25ex,x=.4cm,y=.8cm]
    \foreach \i/\j in {1/2,2/5,3/8,4/3,5/6,6/1,7/4,8/7} {
      \coordinate[label=above:{$\scriptscriptstyle \i$}] (a\i) at (\i,\chordheight);
      \coordinate[label=below:{$\scriptscriptstyle \i$}] (b\i) at (\j,0);
      \draw (a\i) to (b\i);
    }
    \coordinate[label=right:{$\scriptscriptstyle b$}] (b) at (8+\chordheight/2,\chordheight/2);
    \coordinate[label=left:{$\scriptscriptstyle a$}] (a) at (1-\chordheight/2,\chordheight/2);
    \draw[dashed] (a) to (b);
    \draw (1,0) arc (270:90:\chordheight/2);
    \draw (8,\chordheight) arc (90:-90:\chordheight/2);
    \draw (1,0) to (8,0);
    \draw (1,\chordheight) to (8,\chordheight);

    \foreach \i in {1,...,8} {
      \fill (\i,0) circle (1pt); 
      \fill (\i,\chordheight) circle (1pt);
     } 
   \end{tikzpicture}
  \end{center}
 \caption{The chord diagram for $\sigma = [61472583]$.}\label{fig:minorchord}
\end{figure}

Now armed with Theorem~\ref{min}, we describe an infinite family of minimal diagram free knots that are not irreducibly odd.  

\begin{thm}\label{thm:family} Let $\sigma$ be a nonconsecutive permutation on $\{1,2,\cdots,n\}$ (where $n>2$) such that $\sigma_n < n$. Define the permutation
\[
 \hat{\sigma} = [(\sigma_i+1), 2n+3, (\sigma_i+n+2), n+2, (\sigma_i+2n+3), 1]
\]
where $(\sigma_i+x)$ is the sequence $(\sigma_1+x,\sigma_2+x,\cdots,\sigma_n+x)$.  Then $\hat\sigma$ corresponds to a minimal free knot diagram but fails to be irreducibly odd. 
\end{thm}

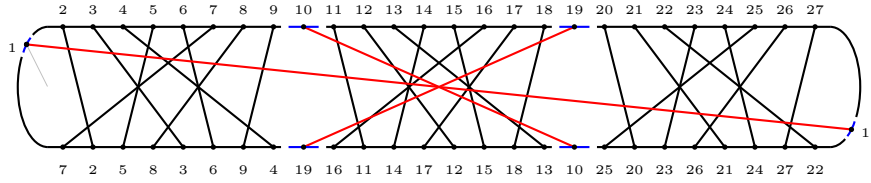
\begin{figure}
\begin{center}
   \newcommand{\chordheight}{2}
   \begin{tikzpicture}[thick,text height=1.5ex,text depth=.25ex,x=.4cm,y=.8cm]
    \draw[blue] (1.5,0) arc (270:90:\chordheight/2);

    \draw[white,ultra thick] (1.5,0) arc (270:145:\chordheight/2);
    \draw[white,ultra thick] (1.5,\chordheight) arc (90:125:\chordheight/2);

    \draw (1.5,0) arc (270:150:\chordheight/2);
    \draw (1.5,\chordheight) arc (90:120:\chordheight/2);

    \draw[blue] (27.5,\chordheight) arc (90:-90:\chordheight/2);

    \draw[white,ultra thick] (27.5,\chordheight) arc (90:-35:\chordheight/2);
    \draw[white,ultra thick] (27.5,0) arc (-90:-55:\chordheight/2);

    \draw (27.5,\chordheight) arc (90:-30:\chordheight/2);
    \draw (27.5,0) arc (-90:-60:\chordheight/2);
   
    \foreach \i/\j in {2/3,3/6,4/9,5/4,6/7,7/2,8/5,9/8,11/12,12/15,13/18,14/13,15/16,16/11,17/14,18/17,20/21,21/24,23/22,22/27,24/25,25/20,26/23,27/26} {
      \coordinate[label=above:{$\scriptscriptstyle \i$}] (a\i) at (\i,\chordheight);
      \coordinate[label=below:{$\scriptscriptstyle \i$}] (b\i) at (\j,0);
      \draw (a\i) to (b\i);
    }
    \foreach \i/\j in {10/19,19/10} {
      \coordinate[label=above:{$\scriptscriptstyle \i$}] (a\i) at (\i,\chordheight);
      \coordinate[label=below:{$\scriptscriptstyle \i$}] (b\i) at (\j,0);
      \draw[red] (a\i) to (b\i);
    }
    \coordinate[label=left:{$\scriptscriptstyle 1$}] (a1) at ($ (1.5,\chordheight/2) + (135:\chordheight/2) $);
    \coordinate[label=right:{$\scriptscriptstyle 1$}] (b1) at ($ (27.5,\chordheight/2) + (315:\chordheight/2) $);
    \draw[red] (a1) to (b1);
    
    \draw[black] (1.5,0) to (9.25,0);
    \draw[blue] (9.5,0) to (10.5,0);
    \draw[black] (10.75,0) to (18.25,0);
    \draw[blue] (18.5,0) to (19.5,0);
    \draw[black] (19.75,0) to (27.5,0);
    
    \draw[black] (1.5,\chordheight) to (9.25,\chordheight);
    \draw[blue] (9.5,\chordheight) to (10.5,\chordheight);
    \draw[black] (10.75,\chordheight) to (18.25,\chordheight);
    \draw[blue] (18.5,\chordheight) to (19.5,\chordheight);
    \draw[black] (19.75,\chordheight) to (27.5,\chordheight);

  \foreach \i in {1,...,27} {
      \fill (a\i) circle (1pt); 
      \fill (b\i) circle (1pt);
     } 
    \fill (1.5,\chordheight/2) -- +(135:\chordheight/2) circle (1pt) [label=left:1];

   \end{tikzpicture}
  \end{center}
 \caption{The construction of $\hat{\sigma}$ from $\sigma = [61472583]$ per Theorem~\ref{thm:family}.}\label{fig:majorchord}
\end{figure}

\begin{proof} We begin by showing that the conditions of Theorem~\ref{min} are satisfied. That is, $\hat\sigma$ is a minimal permutation. 

First, we note that 1 is a crossing element, since $\hat\sigma_{3n+3}=1$ and $\hat\sigma$ is a permutation on $\{1,2,\cdots,3n+3\}$. Furthermore, no other element of $\{1,2,\cdots,3n+3\}$ is a crossing element since $n > 2$.  Also, since $\sigma$ is nonconsecutive, so is $\hat\sigma$. 

To show that (1) holds, we must show that consecutive numbers in the permutation have a difference less than $(3n+3)-2=3n+1$.  We refer to the elements $1, n+2$ and $2n+3$ of $\{1,...,3n+3\}$ as ``red'' and the remaining elements as ``black''. (See Figure~\ref{fig:majorchord} for inspiration.) Let us consider the case where two consecutive numbers are both black and the case where one is red. Note that the case where two consecutive elements are both red doesn't occur, by construction.  Now if both $i-1$ and $i$ are black, then $|\hat{\sigma}_{i-1} - \hat{\sigma}_{i}| \leq n < 3n+1$, as required.  For red elements, note that $\hat{\sigma}_{3n+3} = 1$ and $\hat{\sigma}_{3n+2}=\sigma_n+2n+3 < 3n+3$ (since $\sigma_n< n$) so $|\hat{\sigma}_{3n+2} - \hat{\sigma}_{3n+3}| < 3n-2$.  For the other two red elements, we see that $\hat{\sigma}_{n+2} = 2n+2$ and $\hat{\sigma}_{2n+3} = n+1$ so the longest distance from either to an adjacent element would be at most $\max(|3n+3 - (n+1)|, |2n+2 - 1|) = 2n+2$. %This last sentence is a little confusing.

Now, we turn to (2). Since $\hat\sigma_1=\sigma_1+1$, we have that $2\leq\hat\sigma_1\leq n+1$. Since $\hat\sigma_{3n+2}=\sigma_n+2n+3$ and $\sigma _n<n$, it follows that $2n+4\leq\hat\sigma_{3n+2}\leq 3n+2$. Combining these inequalities, we find that $n+3\leq|\hat\sigma _1-\hat\sigma _{3n+2}|\leq3n$, so $1<|\hat\sigma _1-\hat\sigma _{3n+2}|<3n+1$. Hence, the final minimality property holds.

Since the permutation $\hat\sigma$ is minimal, it corresponds to a minimal free knot diagram. However, the graph associated to $\hat\sigma$ fails to be irreducible. Indeed, the vertices $n+2$ and $2n+3$ are both adjacent to precisely the elements of the set $\{1,n+3,n+4,\cdots,2n+2\}$.  Moreover, if the permutation $\sigma$ that generates $\hat\sigma$ fails to be even, then $\hat\sigma$ will fail to be odd.
\end{proof}

The authors hope that the main result of this paper contributes to the construction of a knot table for free knots. This may be a direction for future work.

%\section{Future Work} %%%Revise%%%

%{\huge I don't think any of this is applicable anymore}
%Our initial goal in studying minimal diagrams of free knots and irreducibly odd graphs was to begin the project of creating a free knot table. While we do not provide such a table here, we have developed several tools that might be used to do so. The next step for this work is to expand on our set of tools to find a complete classification of small crossing number free knots via their minimal diagrams.

%Applications to free knots aside, the graph theoretical concept of reducibility and irreducibility should be investigated further.  Two vertices having the same neighborhood are in a sense indistinguishable.  This can be exploited, for example, in coloring a graph or computing its automorphism group. Much work in this area has yet to be done.

%\vfill
%\pagebreak

%\appendix
%\section{Irreducibly odd graphs on 6 and 8 vertices}

%Table~\ref{table} contains all nonempty irreducibly odd graphs with 8 or fewer vertices.  Chord diagrams are shown for circular graphs, that is, those graphs that are associated to free knots.

%{\huge This table has little value now.  As I see it, there are two options.  We can kill it, or when I get back from vacation, I can delete the graphs which aren't circle graphs, possibly look at irrodd circle graphs on 10 vertices, and maybe we can draw knot diagrams?}

%\begin{table}
  %\begin{center}
    %\input{irrodd_chord}
  %\end{center}
  %\caption{Irreducibly odd graphs on 6 and 8 vertices}\label{table}
%\end{table}

\bibliographystyle{amsplain}
%\nocite{*}
\bibliography{irrodd_apr_17_2014}

\end{document}